\documentclass[a4paper,10pt]{article}

\newfont{\bcb}{msbm10}
\newfont{\matb}{cmbx10}
\newfont{\got}{eufm10}

\usepackage[cp1250]{inputenc}

\usepackage{amsmath,amsthm}
\usepackage{amssymb,latexsym}
\usepackage{enumerate}
\usepackage{latexsym}
      \usepackage[all, knot]{xy}
     \xyoption{arc}

\begin{document}
\title{On continuity of Guo Wuwen function}
\author{Zofia Ambroży, Iwo Biborski}

\footnotetext{2010 MSC. 15A18, 15A29  }

\footnotetext{Key words: matrix theory, nonegative inverse eigenvalue problem, semialgebraic geometry}

\maketitle

\begin{abstract}
We show that the functions $g$ and $g_{s}$ introduced by Guo Wuwen in \cite{[Gu]} are  continuous and semialgebraic. We use this fact to prove that the set $\mathbb{N}_{n}$ of ordered $n$-tuples of real numbers, realizable by nonnegative matrices, is a closed set.  \end{abstract}

\section{Introduction} Let
$$
\mathbb{R}^{n}_{\geq}:=\left\{(\lambda_{1},\dots,\lambda_{n})\in\mathbb{R}^{n}:\lambda_{i}\geq\lambda_{i+1},\,i=1,\dots,n-1\right\}.
$$
For a matrix $A\in\mbox{M}(n,n; \mathbb{R})$, we write $A\geq 0$ if $A$ is a nonnegative matrix, i.e. all entries of $A$ are nonnegative. We denote by $\mbox{Spec}(A)$ the spectrum of $A$.

\vspace{2ex}
\textbf{Definition 1.} We say that $\Lambda\in\mathbb{R}^{n}_{\geq}$ is realizable if there exists a nonnegative matrix $A$ such that $\mbox{Spec}(A)=\Lambda$. We denote the set of all realizable $\Lambda$ by $\mathbb{N}_{n}$. We denote by $\mathbb{S}_{n}$ the set of all $\Lambda$ realizable by a nonnegative symmetric matrix.

\vspace{2ex}
The problem of finding conditions on $\Lambda \in\mathbb{R}^{n}_{\geq}$, such that $\Lambda \in \mathbb{N}_n  \ (\mathbb{S}_n)$  is strictly connected with semialgebraic geometry. Let us recall that semialgebraic subsets of $\mathbb{R}^n$ is the smallest class $\mathcal{SA}_n$ of $\mathbb{R}^n$ such that
\begin{enumerate}[{$(1)$}]
\item if $P \in \mathbb{R}[X_1, \ldots, X_n]$, then
$$\{x \in \mathbb{R}^n \ | \ P(x)=0 \} \in \mathcal{SA}_n \quad {\rm and} \quad
\{x \in \mathbb{R}^n  |  P(x)>0 \} \in \mathcal{SA}_n,$$
\item if $A, B \in \mathcal{SA}_n$, then $A \cap B$, $A \cup B$, $\mathbb{R}^n \backslash A$ are in $\mathcal{SA}_n$.
\end{enumerate}

\noindent For basics of semialgebraic geometry we refer reader to (\cite{[BCR]}, \cite{[Coste]}, \cite{[vdD]}). Now let $A \in M(n,n; \mathbb{R})$ and denote by
$$
f_A(\lambda)= (-1)^n\big( \lambda^n + c_1(A) \lambda^{n-1} + \ldots + c_{n-1}(A)\lambda + c_n(A)\big),
$$
where $c_i: \mathbb{R}^{n^2} \longrightarrow \mathbb{R}$ are suitable polynomials, its characteristic polynomial. Let
$$
a_i: \mathbb{R}^n \ni (\lambda_1, \ldots, \lambda_n)= (-1)^i \sum_{1 \leq j_1 \leq \ldots \leq j_i \leq n}\lambda_{j_1} \cdots \lambda_{j_i} \in \mathbb{R},
$$
for $i=1, \ldots, n$. A subset
$$
\mathbb{M}_n:=\left\{ (\lambda, A) \in \mathbb{R}^n \times \mathbb{R}^{n^2} \ | \ A \geq 0, \ c_i(A)= a_i(\lambda)\right\}
$$
is semialgebraic. Let $\Pi: \mathbb{R}^n \times \mathbb{R}^{n^2} \longrightarrow \mathbb{R}^n$ be a projection on first $n$ coordinates. Of course $\mathbb{N}_n = \Pi(\mathbb{M}_n)$. It follows by Tarski-Seidenberg Theorem (cf. \cite{[BCR]}, Theorem 2.2.1), that the set $\mathbb{N}_n$ is semialgebraic. In a similar way we can show that the set $\mathbb{S}_n$ is semialgebraic.

\vspace{2ex}
\textbf{Definition 2(Guo Wuwen, \cite{[Gu]}).}
We define the following function $g:\mathbb{R}^{n-1}_{*}\rightarrow\mathbb{R}$:
\begin{gather*}
g(\lambda_{2},\dots,\lambda_{n})=\inf_{t\in\mathbb{R}}\left\{t\geq\max_{i=1, \ldots, n-1}|\lambda_{i}|,\, \forall_{\delta\geq0}(t+\delta,\lambda_{2},\dots,\lambda_{n})\in\mathbb{N}_{n}\right\}.
\end{gather*}
We denote by $g_{s}$  the analogical function for symmetric realization. Guo Wuwen proved, that
$$
\forall_{t \in \mathbb{R}} \quad t > g(\lambda_2, \ldots, \lambda_n) \ \Rightarrow \ (t, \lambda_2, \ldots, \lambda_n) \in \mathbb{N}_n,
$$
and similar fact for the function $g_s$.

\vspace{2ex}
We will now show that $g$ is a semialgebraic function.  Let us recall that a mapping $f:A \longrightarrow B$, where $A \subset \mathbb{R}^n$ and $B \subset \mathbb{R}^m$ are semialgebraic, is said to be semialgebraic, if its graph
$$
graph(f)= \left\{ (x,y) \in A \times B \ | \ y=f(x)\right\}
$$
is a semialgebraic subset of $\mathbb{R}^n \times \mathbb{R}^m$.
We define
$$
B:=\left\{ (\lambda, \epsilon, \lambda_1, \ldots, \lambda_n) \in \mathbb{R}^{n+2} \ | \ (\lambda, \lambda_2, \ldots, \lambda_n) \in \mathbb{N}_n, \
\epsilon > 0, \ \lambda_1=\lambda + \epsilon \right\}.
$$
For every $(\lambda_2, \ldots, \lambda_n) \in \mathbb{R}^n_{\geq}$ we have
$$
\big(\widetilde{\Pi}(B)\backslash \mathbb{N}_n\big) \cap \big(\mathbb{R} \times (\lambda_2, \ldots, \lambda_n)\big) = \big(g(\lambda_2, \ldots, \lambda_n), \lambda_2, \ldots, \lambda_n \big),
$$
where $\widetilde{\Pi}: \mathbb{R} \times \mathbb{R} \times \mathbb{R}^{n} \ni (x, y ,z) \longrightarrow z \in \mathbb{R}^n$. We get that
$graph(g)= \widetilde{\Pi}(B)\backslash \mathbb{N}_n$.

\vspace{2ex}
Our purpose is to prove that $g$ and $g_{s}$ are continuous semialgebraic functions. Then we use this fact to prove that the set of all realizable real spectra is closed. We also consider the following problem: if the spectrum $\Lambda$ is realizable and all eigenvalues are pairwise different then $\Lambda$ is symmetrically realizable? The answer is negative. The solution comes from continuity of $g$ and $g_{s}$ supported with a~ proper example.

\section{Main results.}
Let us recall the following

\vspace{2ex}
\textbf{Corollary 1(Guo Wuwen, \cite{[Gu]}, Corollary 3.2).}\textit{ Let $n \in \mathbb{N}$  and assume that
\begin{enumerate}[{$(1)$}]
\item $(\lambda_1; \lambda_2, \ldots, \lambda_n) \in \mathbb{N}_n$,
\item $\epsilon_i \in \mathbb{R}$ for $i=1, \ldots, n$, and $\epsilon_1= \sum_{i=2}^{n}|\epsilon_i|$.
\end{enumerate}
Then
$$
(\lambda_1 + \epsilon_1; \lambda_2 + \epsilon_2, \ldots, \lambda_n + \epsilon_n) \in \mathbb{N}_n.
$$}

We use Corollary 1 to prove the following

\vspace{2ex}
\textbf{Theorem 2.}\textit{ The function $g$ is continuous.}
\begin{proof}
We shall show that $g$ is upper continuous and lower continuous.

\vspace{2ex}
\textit{Case 1.} $g$ is upper continuous. Let $(\lambda_2, \ldots, \lambda_n) \in \mathbb{R}^{n-1}$ and let $\epsilon >0$. Take $(\tilde{\lambda}_2, \ldots, \tilde{\lambda}_n)$ such that
$$
\tilde{\epsilon}=\sum_{i=2}^{n} |\lambda_i - \tilde{\lambda}_i|< \epsilon.
$$
By corollary,
$$
g(\tilde{\lambda}_2, \ldots, \tilde{\lambda}_n) \leq \lambda_1 + \tilde{\epsilon} \leq \lambda_1 + \epsilon.
$$
It follows, that
$$
\limsup g(\tilde{\lambda}_2, \ldots, \tilde{\lambda}_n) \leq g(\lambda_2, \ldots, \lambda_n).
$$

\vspace{2ex}
\textit{Case 2.} $g$ is lower continuous. Let $\lambda=(\lambda_2, \ldots, \lambda_n) \in \mathbb{R}^{n-1}$ and suppose, that
$$
\sigma_0 =\liminf_{\delta \rightarrow \lambda} g(\tilde{\lambda}) < g(\lambda).
$$
There exists $\delta^k \in \mathbb{R}^{n-1}, k \in \mathbb{N}$ such that
$$
\sum_{i=2}^{n}|\delta_{i}^{k} - \lambda_{i}|\leq \frac{1}{n}, \ {\rm and} \ \ g(\delta^{k}) \underset{k \rightarrow \infty}{\longrightarrow} \sigma_0.
$$
Let $\widetilde{\sigma}= g(\lambda) - \sigma_0$. Take such $k_0 \in \mathbb{N}$, that for every $k \geq k_0$
$$
\sum_{i=2}^{n}|\delta_{i}^{k} - \lambda_{i}|\leq \frac{\widetilde{\sigma}}{2}, \quad  g(\delta^{k}) < g(\lambda) - \frac{\widetilde{\sigma}}{2} \quad {\rm and} \quad g(\delta^{k}) + \frac{\widetilde{\sigma}}{2} \geq \max\{|\lambda_i|\}.
$$
Let $\epsilon_i \in \mathbb{R}, i=2, \ldots, n$ be such that $\lambda_i - \delta_{i}^{k_0}= \epsilon_i$. Then
$$
(g(\delta^{k_0}), \delta^{k_0}_{2}, \ldots, \delta^{k_0}_{n}) \in \mathbb{N}_n,
$$
and by corollary
$$
(g(\delta^{k_0})+ \frac{\widetilde{\sigma}}{2}, \lambda_2, \ldots, \lambda_n) \in \mathbb{N}_n.
$$
But
$$
g(\delta^{k_0}) + \frac{\widetilde{\sigma}}{2} < g(\lambda) - \frac{\widetilde{\sigma}}{2} + \frac{\widetilde{\sigma}}{2} = g(\lambda),
$$
a contradiction.
\end{proof}

We use a different argument for the function $g_{s}$. The reasoning for symmetric realizability is based on a standard fact that for any symmetric matrix $A$ there exist a diagonal matrix $D$ and an orthogonal matrix $U\in O(n)$ such that $U^{T}DU=A$, where $X^{T}$ denotes the transpose of a given matrix $X$ and $O(n)$ is the orthogonal group of matrices of dimension $n$. We have the following

\vspace{2ex}
\textbf{Theorem 3.}\textit{ The function $g_{s}:\mathbb{R}^{n-1}\rightarrow\mathbb{R}$ is continuous.}
\begin{proof}Let $\Lambda=(\lambda_{2},\dots,\lambda_{n})$  and $\Lambda_{k}:=(\lambda_{k, 2},\dots,\lambda_{k, n})$ be the sequence of $(n-1)$-tuples such that $\lim_{k\rightarrow\infty}\Lambda_{k}=\Lambda$. Let $t=g_{s}(\Lambda)$ and $t_{k}=g_{s}(\Lambda_{k})$. We can divide $\Lambda_{k}$ into two subsequences $\Lambda_{k_{1}}$ and $\Lambda_{k_{2}}$  such that $g_{s}(\Lambda_{k_{1}})> g_{s}(\Lambda)$ and $g_{s}(\Lambda_{k_{2}})\leq g_{s}(\Lambda)$. Thus we can just assume that $t_{k}>t$ for all $k\in\mathbb{N}$ or that $t_{k}\leq t$ for all $k\in\mathbb{N}$.

\vspace{2ex}
\textit{Case 1.} Suppose that $t_{k}>t$ for all $k\in\mathbb{N}$. Let $\varepsilon>0$. Thus there exists $U\in O(n)$ such that $A:=U^{T}\mbox{diag}(t+\varepsilon,\lambda_{2},\dots,\lambda_{n})U$ is a positive matrix. Consider a sequence of matrices $A_{k}:=U^{T}\mbox{diag}(t+\varepsilon,\lambda_{k, 2},\dots\lambda_{k, n})U$. It is clear that $\lim_{k\rightarrow\infty}A_{k}=A$, thus there exists $K\in\mathbb{N}$ such that, for all $k>K$, $A_{k}$ is a positive matrix. Whence, for $k>K$, $t_{k}\leq t+\varepsilon$. Finally, for all $\varepsilon>0$, there exists $K\in\mathbb{K}$ such that $t_{k}\in(t,t+\epsilon]$ and thus $t_{k}\rightarrow t$ while $k\rightarrow\infty$.

\vspace{2ex}
\textit{Case 2.} Suppose that $t_{k}\leq t$ for all $k\in\mathbb{N}$. Thus $\{t_{k}\}\subset[0,t]$. There exists a sequence $U_{k}$ of orthogonal matrices such that
\begin{gather*}
U_{k}^{T}\mbox{diag}(t_{k},\lambda_{k 2},\dots,\lambda_{k n})U_{k}
\end{gather*}
is a nonnegative matrix. Since $\{t_{k}\}\subset[0,t]$, each subsequence of $\{t_{k}\}$ has an accumulation point. Suppose that $\tilde{t}$ is an accumulation point of $\{t_{k}\}$, and $\{t_{k_{s}}\}$ converges to $\tilde{t}$. Since $O(n)$ is a compact set, thus there exists a convergent subsequence $\{U_{k_{\nu}}\}$ of $\{U_{k_{s}}\}$. Let

\begin{gather*}
\lim_{k_{\nu}\rightarrow\infty}U_{k_{\nu}}:=U.
\end{gather*}

It is clear that $\{t_{k_{\nu}}\}$ converges to $\tilde{t}$. Thus

\begin{gather*}
A_{k_{\nu}}:=U^{T}_{k_{\nu}}\mbox{diag}(t_{k_{\nu}},\lambda_{k_{\nu}, 2},\dots\lambda_{k_{\nu}, n})U_{k_{2}}
\end{gather*}
form a sequence of nonnegative matrices with limit $A:=U^{T}\mbox{diag}(\tilde{t},\lambda_{2},\dots,\lambda_{n})U$. Obviously $A$ is nonnegative as the limit of the sequence of nonnegative matrices, and thus $\tilde{t}\geq t$. Since $\{t_{k}\}\subset[0,t]$, $\tilde{t}\leq t$ and therefore $\tilde{t}=t$. Finally, each accumulation point of $\{t_{k}\}$ is equal to $t$ and thus $\lim_{k\rightarrow\infty}t_{k}=t$.

This ends the proof that $g_{s}$ is a continuous function.
\end{proof}

\textbf{Corollary 4.}\textit{ The sets $\mathbb{N}_{n}$ and $\mathbb{S}_{n}$ are closed semialgebraic sets.}
\begin{proof}
It is clear that $\mathbb{N}_{n}$ and $\mathbb{S}_{n}$ are semialgebraic.

Let $\Lambda_{k}:=(\lambda_{k, 1 },\dots,\lambda_{k, n})\in\mathbb{N}_{n}$ be the sequence convergent to $\Lambda:=(\lambda_{1},\dots,\lambda_{n})$. Let $\Lambda_{k}':=(\lambda_{k, 2},\dots,\lambda_{k, n})$ and $\Lambda':=(\lambda_{2},\dots,\lambda_{n})$. Thus
\begin{gather*}
 \lambda_{k, 1}\geq g(\Lambda_{k}').
\end{gather*}
Since $\lambda_{k, 1}$ is convergent to $\lambda_{1}$ and $g$ is continuous, thus
\begin{gather*}
 \lambda_{1}\geq g(\Lambda'),
\end{gather*}
and therefore $\Lambda\in\mathbb{N}_{n}$. Finally $\mathbb{N}_{n}$ is closed. For $\mathbb{S}_{n}$ the proof is similar.
\end{proof}

\vspace{2ex}
Let $\gamma=(\gamma_1, \ldots, \gamma_n): (0, \epsilon) \longrightarrow \mathbb{N}_n$ (where $\epsilon>0$) be a semialgebraic continuous curve, such that $\gamma$ is bounded. Then $\gamma(0):= \lim_{t \rightarrow 0^+} \gamma(t)$ exists and $\gamma(0) \in \mathbb{N}_n$ by Corollary 4 and the fact, that every semialgebraic curve $\alpha: (0, \delta) \longrightarrow \mathbb{R}^m$ (where $\delta >0, m \in \mathbb{N}$) has $\lim_{s \rightarrow 0^+} \alpha(s)$. We have the following

\vspace{2ex}
\textbf{Proposition 5.}\textit{Let $\gamma=(\gamma_1, \ldots, \gamma_n): (0, \epsilon) \longrightarrow \mathbb{N}_n$ be a bounded, continuous semialgebraic curve. Then there exists $\widetilde{\gamma}=(\widetilde{\gamma}_1, \ldots, \widetilde{\gamma}_{n^2 +n}):(0, \epsilon') \longrightarrow \mathbb{M}_n$, where $\epsilon' \leq \epsilon$, such that
\begin{enumerate}[{$(1)$}]
\item $\widetilde{\gamma}_1= \gamma_1|_{(0, \epsilon')}, \ldots, \widetilde{\gamma}_n= \gamma_n|_{(0, \epsilon')},$
\item $\widetilde{\gamma}$ is continuous,
\item $\widetilde{\gamma}(0):=\lim_{t \rightarrow 0^+} \widetilde{\gamma}$ exists and $\big(\widetilde{\gamma}_{n+1}(0), \ldots, \widetilde{\gamma}_{n+n^2}(0)\big) \in \mathbb{M}_{n, \gamma(0)}$.
\end{enumerate}}

\begin{proof}
Let us recall that if $\Lambda \in \mathbb{N}_n$, then there exists $M_{\Lambda} \in \mathbb{M}_{n, \Lambda}$ such that $e=(1, \ldots, 1) \in \mathbb{R}^n$ is an eigenvector for $\rho(\Lambda)=\max_{\lambda \in \Lambda}\{|\lambda|\} $ (cf. \cite{[Gu]}, Lemma 2.2). In particular, every coefficient of $M_{\Lambda}$ is not greater than $\rho(\Lambda)$. Let us define
$$
\mathbb{M}^{\rho}_{n}:=\left\{ (\Lambda, A) \in \mathbb{R}^{n} \times \mathbb{R}^{n^2} \ | \
(\Lambda, A) \in \mathbb{M}_n, \ A e = \rho(\Lambda)\right\}.
$$
The set $\mathbb{M}^{\rho}_{n}$ is semialgebraic. By the Definable Choice Theorem (\cite{[vdD]}, Chapter 6, 1.1.2) there exists
$$
\widetilde{\gamma}:(0, \epsilon) \longrightarrow \mathbb{M}_n, \quad {\rm such \ that} \quad \widetilde{\gamma}_1= \gamma_1, \ldots, \widetilde{\gamma}_n= \gamma_n.
$$
Possibly taking $0 < \epsilon' \leq \epsilon$ we may assume, that $\widetilde{\gamma}|_{(0, \epsilon')}$ is continuous. Since $\gamma$ is bounded, $\widetilde{\gamma}$ is also bounded and there exists $\widetilde{\gamma}(0):=\lim_{t \rightarrow 0^+} \widetilde{\gamma}$. By Polynomial Roots Continuity Theorem,
$\big(\widetilde{\gamma}_{n+1}(0), \ldots, \widetilde{\gamma}_{n+n^2}(0)\big) \in \mathbb{M}_{n, \gamma(0)}$.
\end{proof}

\vspace{2ex}
It is known that the real and symmetric nonnegative inverse eigenvalue problems are different for the dimension equal or higher than 5(\cite{[JLL]}). One can ask if a given realizable spectrum $\Lambda:=(\Lambda_{1},\dots,\Lambda_{n})$ such that $\lambda_{i}\neq\lambda_{j}$ for $i\neq j$, is symmetrically realizable. The continuity of $g$ and $g_{s}$ can be used to prove that this assertion is false.

\vspace{2ex}
\textit{Example 1.}
Consider the well known example od dimension $5$(\cite{[LM]}, \cite{[HL]}, \cite{[Me]}). Put $\Lambda=(3,3,-2,-2,-2)$. It is known from Perron-Frobenius theory that $\Lambda$ is not realizable, but Loewy and Hartwig proved that $(4,3,-2,-2,-2)$ is symmetrically realisable. On the other hand Meehan proved that there exists $t\in(0,1)$ such that $(3+t,3,-2,-2,-2)$ is realizable(\cite{[Me]}). Thus $g(\Lambda')\in(3,4)$ and $g_{s}(\Lambda')=4$, where $\Lambda'=(3,-2,-2,-2)$. Let
\begin{gather*}
\Lambda'_{k}:=(3+1/k,-2+1/k,-2-2/k,-2+3/k),\,\,\,k\in\mathbb{N}.
\end{gather*}
It is clear that $\Lambda'_{k}$ converges to $\Lambda'$. By Theorem 2 and Theorem 3, $g(\Lambda'_{k})<4$ for sufficiently large $k$, and $\lim_{k\rightarrow\infty}g_{s}(\Lambda'_{k})=4$. Therefore, by Corollary 1, there exists $k\in\mathbb{N}$ such that $3+t+7/k<4$, $\sigma:=(3+t+7/k,3+1/k,-2+1/k,-2+2/k,-2+3/k)\in\sigma\in\mathbb{N}_{5}$, $\sigma\notin\mathbb{S}_{5}$ and the elements of $\sigma$ are necessarily pairwise different.


\begin{thebibliography}{99}
\bibitem{[BCR]} J. Bochnak, M. Coste, M.-F. Roy. \textit{Real algebraic geometry}. Springer, 1998.
\bibitem{[Coste]} M. Coste. \textit{An introduction to semialgebraic geometry}. Dottorato di Ricerca in Matematica, Dip. Univ. Pisa. Instituti Editoriali e Poligrafici Internazionali, 2000.
\bibitem{[vdD]} L. van den Dries. \textit{Tame topology and o-minimal structures}. London Math. Soc. Lecture Note 248. Cambridge Univ. Press, 1998.
\bibitem{[Gu]} Guo Wuwen. \textit{Eigenvalues of nonnegative matrices}. Linear Algebra Appl. \textbf{266} (1997), 261-270.
\bibitem{[HL]} R. Hartwig and R. Loewy. Unpublished (1989)
\bibitem{[JLL]} C.R. Johnson, T.J. Laffey, R. Loewy. \textit{The real and the symmetric
nonnegative inverse eigenvalue problems are different}. Proc. of the AMS \textbf{124} (1996), 3647-3651.
\bibitem{[LM]} R.Loewy, J.J. McDonald. \textit{The Symmetric Nonnegative Inverse
Eigenvalue Problem for $5\times5$ Matrices}.  Linear Algebra Appl. \textbf{393} (2004), 275-298
\bibitem{[Me]} M.E. Meehan. \textit{Some results on matrix spectra}. PhD Thesis, University
College Dublin, 1998.


\end{thebibliography}
\end{document}